\documentclass[11pt]{article}

\usepackage{amssymb}
\usepackage{amsmath}

\newcommand{\mapto}{\to}
\newcommand{\R}{\mathbf R}
\newcommand{\Z}{\mathbf Z}
\newcommand{\Q}{\mathbf Q}

\newcommand{\I}{\mbox{\cal I}}
\newtheorem{Theorem}{Theorem}

\newtheorem{Lemma}[Theorem]{Lemma}

\newtheorem{Corollary}[Theorem]{Corollary}

\newcommand{\sgn}{\mbox{sgn}}

\begin{document}

\begin{center}
{\bf \Large Pfaffian formulae for one dimensional coalescing and annihilating systems}
\vskip 12pt

Roger Tribe and Oleg Zaboronski

Department of Mathematics, University of Warwick,
CV4 7AL, UK. r.p.tribe@warwick.ac.uk, olegz@maths.warwick.ac.uk



Submitted September 22, 2010, accepted September 16, 2011.

\end{center}

\begin{abstract}
The paper considers instantly coalescing, or instantly annihilating, systems of
one-dimensional Brownian particles on the real line.
Under maximal entrance laws, the distribution of the particles at a fixed time is
shown to be Pfaffian point processes closely related to the Pfaffian point process describing one dimensional distribution of
real eigenvalues in the real Ginibre ensemble of random matrices. As an application, an exact large time asymptotic for the
$n$-point density function for coalescing particles is derived.
\end{abstract}

\section{Introduction and summary of main results} \label{s1}
The study of single species reaction diffusion systems $A+A \to A$ (coalescence)
and $A + A \to 0$ (annihilation) originated in non-equilibrium statistical mechanics (see \cite{privman}), but
has now a large mathematical literature (see, for example, \cite{arratia}, \cite{bramson+griffeath}, \cite{bramson+lebowitz},
\cite{kesten+vandenberg}). In one dimension the systems exhibit strongly
non-mean field behaviour due to correlation effects. In this paper
we give several examples showing that this correlation structure can be encoded algebraically
in a Pfaffian structure. Note that the embedding of annihilating random walks
as domain boundaries for a Glauber model makes Pfaffian formulae quite reasonable due
to the free fermion structure of the Glauber model (see Felderhof \cite{felderhof}).

We examine the asymptotics for the $n$-particle
density function $ \rho_t^{(n)}$ for (instantly) coalescing Brownian motions on $\R$ defined by
\[
P \left[ \mbox{there exist particles in $dx_1,\ldots,dx_n$ at time $t$} \right]
= \rho_t^{(n)}(x_1,\ldots,x_n) \, dx_1 \ldots dx_n.
\]
In \cite{mrtz} we showed, for $n \geq 1$, $t_0,L>0$ and for a variety of initial conditions, the bounds
\begin{equation} \label{oldbounds}
0<c_1(n,L,t_0) \leq \frac{ \rho^{(n)}_t(x_1,\ldots,x_n)}{t^{-\frac{n}{2} - \frac{n(n-1)}{4}}
\prod_{1 \leq i < j \leq n} |x_i-x_j|}
 \leq c_2(n,L,t_0) < \infty,
\end{equation}
for $t \geq t_0$ and $|x_i| \leq L t^{1/2}$, where the constant $c_1$ will depend also
on the initial condition. The non-linear factor in the power of $t$ illustrates the non mean-field behaviour
due to correlations.

In this paper we show that, under the maximal entrance law, the true asymptotic holds in (\ref{oldbounds})
as $t \to \infty$, and identify the limiting constant as the Pfaffian of a certain matrix.
The maximal entrance law corresponds intuitively to starting with every point occupied, and can be constructed
as the limit of initial Poisson distributions with increasing intensities. This initial condition is natural since,
as explained in section \ref{s2.3},  started from a large class of other initial
conditions the distributions at time $t$ become close, as $t \to \infty$, to those of the maximal entrance law.
\begin{Theorem} \label{T1}
Under the maximal entrance law for coalescing Brownian motions,
\[
\sup_{|x_{i}| << t^{1/2}} \left|
\frac{\rho^{(2n)}_t(x_1,\ldots,x_{2n})}
{t^{-n-\frac{n(2n-1)}{2}} \prod_{1 \leq i < j \leq 2n} |x_i-x_j|}
- (4\pi)^{-n/2} \mbox{Pf}\,(J^{(2n)}(\phi)) \right| \to 0 \quad \mbox{as $t \to \infty$,}
\]
where $|x_{i}| <<  t^{1/2}$ means that we may take the supremum over any positions
$(x_{i}(t))$ provided that $\sup_i |x_i(t)|t^{-1/2} \to 0$ as $t \to \infty$, and where
$ \mbox{Pf}\,(J^{(2n)}(\phi))$ is the
Pfaffian of the $2n \times 2n$ anti-symmetric matrix $J^{(2n)}(\phi)$ with entries
\[
J^{(2n)}_{ij}(\phi) = (-1)^{j-1} \frac{1}{(i-1)! (j-1)!}
\frac{d^{i+j-2} \phi}{d x^{i+j-2}} (0) \quad \mbox{for $1 \leq i <j \leq 2n$,}
\]
where $\phi(z) = z \exp(-z^2/4)$. Under the analogous maximal entrance law for annihilating
Brownian motions, the same limit holds with $(4\pi)^{-n/2}$ replaced by $(64 \pi)^{-n/2}$.
\end{Theorem}
The presence of the Pfaffian in this asymptotic is a reflection that under the maximal
entrance law the particle positions, at a fixed time, form a Pfaffian point process
(see the start of section \ref{s3} for a definition).
% Under the maximal entrance law there is a scaling relation: if we write $N_t$ for the empirical measure
% generated by the particle positions at time $t$, then  $N_{c^2t}$ has the same distribution
%as the stretched image of the measure $N_t$ under the map $x \to cx$.
%
%
\begin{Theorem} \label{T2}
Under the maximal entrance law for coalescing Brownian motions, the particle
positions at time $t$ form a Pfaffian point process with kernel
$t^{-1/2} K(xt^{-1/2},yt^{-1/2})$, where
\[
K(x,y) = \left( \begin{array}{cc}
-F''(y-x) & -F'(y-x) \\
F'(y-x) & \sgn(y-x) F(|y-x|)
\end{array} \right)
\]
and  $F(x) = \pi^{-1/2} \int^{\infty}_x e^{-z^2/4} dz$. (Here $\sgn(z) = 1$ for $z>0$, $\sgn(z) = -1$ for $z<0$ and $\sgn(0) =0$.)
Under the analogous maximal entrance law for annihilating Brownian motions, the particle
positions at time $t$ form a Pfaffian point process with kernel $\frac12 t^{-1/2} K(xt^{-1/2},yt^{-1/2})$.
\end{Theorem}
The annihilating versions of Theorems \ref{T1} and \ref{T2}, that is for (instantly) annihilating particles, can be deduced from
the thinning relation that connects coalescing and annihilating systems
(see section \ref{s2.1}).

Many probabilities for the fixed $t$ distributions are given by formulae using
Pfaffians, and there are many places to start when proving these formulae. We choose to start
by considering the following basic fact for product moments for annihilating systems, from which we will deduce all
the other Pfaffian fromulae.
\begin{Theorem} \label{T3}
Consider the product moments for annihilating Brownian motions, defined by
\[
m^{(n)}_t(x_1,\ldots,x_{n}) =  E^{A}_{(x_1,\ldots,x_{n})} \left[ \prod_{i\in I_t} g(X^i_t) \right],
\]
for bounded measurable $g$, where $(x_1, \ldots,x_{n})$ lists the initial positions of the
annihilating Brownian motions on $\R$,
and $(X^i_t: i \in I_t)$ list the positions of any particles that remain at time $t$
(and an empty product is taken to have value $1$).
Then for $x_1 < x_2 < \ldots < x_{2n}$, the even moments $m^{(2n)}_t(x)$ are given by
\begin{equation} \label{evenpf}
m^{(2n)}_t(x_1,\ldots,x_{2n}) = \mbox{Pf} \, \left(m^{(2)}_t(x_i,x_j): 1 \leq i < j \leq 2n \right)
\end{equation}
where the right hand side is the Pfaffian of the $2n \times 2n$ anti-symmetric matrix with
entries $m^{(2)}_t(x_i,x_j)$ above the diagonal.
\end{Theorem}
Note these Pfaffians are in variables that determine the initial conditions, allowing   us to use
p.d.e. methods to characterize these moments. Indeed, the product moments satisfy a closed
system of heat equations (with suitable boundary
conditions), and we will verify Theorem \ref{T3} by simply checking that the Pfaffian
uniquely satisfies this system. Markov time-reversal duality (see section \ref{s2.2}) then immediately
implies that certain empty interval formulae
\[
P \left[ \mbox{the intervals $(a_1,a_2)$, $(a_3,a_4)$,...,$(a_{2m-1},a_{2m})$ are empty at time $t$} \right]
\]
for coalescing systems are given by a
Pfaffian, where the Pfaffian is now in the variables $a_1<a_2< \ldots <a_{2m}$
that determine the end points of the target intervals.
This quickly leads to the identification of the Pfaffian point process kernel $K(x,y)$.

We concentrate on Brownian particles but, as we note later, we expect many of our Pfaffian formulae to
hold for a large variety of spatial motion processes, and the Pfaffian structure seems to
arise from two basic underlying mechanisms: linearly ordered particle motion and instantaneous
reactions.
\subsection{Relation between annihilating Brownian motions and the real Ginibre ensemble of random matrices.}
The Pfaffian point process defined in Theorem \ref{T2} has been originally discovered in the context of
random matrices.\footnote{We are grateful to the anonymous referee for bringing this fact to our attention.}
Namely, consider real Ginibre ensemble \cite{Ginibre} defined by the following probability measure on
the space of real $N\times N$ matrices:
\begin{eqnarray}
\mu(d\mathbf{M})=\frac{1}{(2\pi)^{N^2/2}}e^{-\frac{1}{2}Tr(\mathbf{M}^T\mathbf{M})}\lambda_{N\times N}(d\mathbf{M}),
\end{eqnarray}
where $\lambda_{N\times N}$ is Lebesgue measure on $\R^{N\times N}$.
Even though the real Ginibre ensemble is a classical  matrix model, the eigenvalue correlation
functions have been computed only recently, see \cite{borodin_sinclair}, \cite{sinclair}, \cite{nagao},
\cite{sommers} and \cite{sommers1}.

It turns out that the pfaffian point process corresponding to one-dimensional distributions of annihilating Brownian motions
is equivalent to the pfaffian point process describing the law of $real$ eigenvalues of Ginibre in the limit $N\rightarrow \infty$.
Namely, comparing the statement of Theorem \ref{T2} with Corollary 9 of \cite{borodin_sinclair} we arrive at the following
conclusion:
\begin{Corollary}\label{C1}
The one-dimensional law of particle positions for the system of annihilating Brownian motions on $\mathbb{R}$ at time $t>0$ under the maximal entrance law is a Pfaffian point process with
the kernel
\begin{eqnarray}
K_t^{ABM}(x,y)=\frac{1}{\sqrt{2t}}K^{Ginibre}_{rr}\left(\frac{x}{\sqrt{2t}}, \frac{y}{\sqrt{2t}}\right),
\end{eqnarray}
where $K^{Ginibre}_{rr}$ is the $N\rightarrow \infty$ limit of the kernel of the Pfaffian point process characterizing the
law of real eigenvalues in the real Ginibre ensemble.
\end{Corollary}
In other words, the one-dimensional law of annihilating Brownian
motions under the maximal entrance law initial conditions is equivalent to the $N=\infty$ limiting
law of real eigenvalues of a real matrix with independent normal entries.

Corollary \ref{C1} suggests that real eigenvalues of real matrix-valued Brownian motion might behave like a system of
one-dimensional annihilating Brownian motions.  In fact, numerical evidence we accumulated up to date compels us to make the following conjecture.
\\
\\
{\bf Conjecture.} {\it Under the maximal entrance law, all finite-dimensional distributions of particle
 positions for a system of annihilating Brownian motions on $\R$
coincide with $N \rightarrow \infty$ limit of multi-time correlation functions of real eigenvalues of
$gl_{\mathbf{R}}(N)$-valued Brownian motion.}\\
\\
Here $gl_{\mathbf{R}}(N)$ denotes the linear space of all $N\times N$ matrices with
real entries.

\section{Brief review of some facts for one-dimensional
coalescing and annihilating Brownian motions}  \label{s2}

We consider, at first, initial conditions that have only finitely many particles.
This paper describes only the one dimensional time distributions, that is at a fixed $t>0$, of
any remaining particles. We list the positions of the particles at time $t$ as $(X^i_t : i \in I_t)$.
The exact details of the labeling system $I_t$ will not be important for us, and indeed our results all
relate only to the empirical measure $N_t$ defined by
\[
N_t(A) = \sum_{i \in I_t} \chi(X^i_t \in A) \quad \mbox{for measurable $A \subseteq \R$.}
\]
For the case of annihilating particles, if the initial number of particles is even then it remains so
for all time. To a list $(x_i)$ of an even number $2n$ of disjoint positions we may associate the open set
\[
S((x_i)) = (\hat{x}_1, \hat{x}_2) \cup \ldots \cup  (\hat{x}_{2n-1},\hat{x}_{2n})
\]
where $\hat{x}_1 < \ldots < \hat{x}_{2n}$ are the ordered positions.
Some of the formulae for annihilating particles
are then most easily stated in terms of the set valued process
\[
S_t = S \left((X^i_t:i \in I_t) \right).
\]
\noindent
\textbf{Notation. }
We write $P^{C}_{(x_1,\ldots,x_n)}$ to indicate that we are considering (instantly)
coalescing Brownian motions started from initial positions $x_1,\ldots,x_n$. When the particles are
annihilating we change the superscript from $C$ to $A$. When the initial positions are random we change the subscript to $\Xi$,
where $\Xi$ is the law of $(X^i_0:i \in I_0)$.
\subsection{The thinning relation} \label{s2.1}

The formulae about coalescing systems in the paper will always come with an analogue for annihilating systems.
The close link between the two systems has often been observed. For this paper the formulae can usually be
derived from the following thinning relation. For a list of positions $(x_1,\ldots,x_n)$ we let
$\Theta(x_1,\ldots,x_n)$ be the random subset
formed by thinning at rate $1/2$, that is by removing each position independently with probability $1/2$. We may also
thin a random set of positions, for example $\Theta(X^i_t:i \in I_t)$, with the understanding that the randomness in the thinning is
independent of the randomness in the set of positions. We write $\Theta(\Xi)$ for the law of the
thinned set of positions that initially have law $\Xi$. Then the thinning relation between
coalescing and annihilating Brownian motions is the following equality in distribution:
\begin{equation} \label{thinning}
\mbox{$(X^i_t:i \in I_t)$ under $P^{A}_{\Theta(\Xi)}$} \quad \stackrel{\mathcal{D}}{=}
\quad \mbox{$\Theta(X^i_t:i \in I_t)$ under $P^{C}_{\Xi}$} .
\end{equation}
Such a thinning relation is discussed in Arratia \cite{arratia} for the scaled limit of
reacting random walks, and is related to results in many later papers. There is a
simple colouring proof (see ben Avraham and Brunet \cite{avraham2})
that bears repetition here. After the paths of a coalescing system have been realized, independently add random
colours as follows. Initially colour each particle red or blue independently with probability
$1/2$. At coalescences the colours evolve according to the rules $R+R \mapto R$, $B+B \mapto R$ and
$R+B \mapto B$. Then the resulting system of blue particles evolves as an annihilating system.
Moreover the colour of a particle at time $t$ depends on whether there were a odd or even
number of ancestors at time zero that were coloured blue. Since distinct particles have disjoint sets of
ancestors, the colour of all particles at any time $t>0$ remains independently red or blue
with equal probability. The thinning relation follows. This argument makes it clear that the result
holds much more widely, since the exact nature of the motion process is not relevant, nor is the
mechanism of reaction (for example it holds for delayed reactions, when the reactions are
controlled by the intersection local time).
\subsection{Duality formulae} \label{s2.2}

We use two duality formulae.
For $a_1<a_2< \ldots <a_{2m}$ let $I_k = (a_k,a_{k+1})$ for $k=1,\ldots,2m-1$. Then
for disjoint $(x_i)$
\begin{equation} \label{cbmduality}
P^{C}_{(x_i)} \left[ N_t(I_1)=N_t(I_3) = \ldots =N_t(I_{2m-1})=0 \right]
=
P^{A}_{(a_i)} \left[ S_t \cap (x_i) = \emptyset \right].
\end{equation}
The annihilating analogue of this is, writing  $|A|$ for the cardinality of a set $A$,
\begin{equation} \label{abmduality}
E^{A}_{(x_i)} \left[ (-1)^{N_t(I_1 \cup I_3 \cup \ldots \cup I_{2m-1})} \right]
=
E^{A}_{(a_i)} \left[ (-1)^{| S_t \cap (x_i)|} \right].
\end{equation}
There are various ways to see these formulae, but for coalescing systems
a key construction is the Brownian web and its coupling with the dual Brownian web, first considered
by Arratia and explored in Toth and Werner \cite{toth+werner} (and subsequent papers).
We need only part of the Brownian web as follows. For a fixed $t>0$, there is a system of
coalescing Brownian motions starting from every rational $x$ and running over the time interval $[0,t]$,
and a coupled system of backwards coalescing Brownian paths starting at time $t$ at all $x \in \Q$ and
running back to time zero. In fact, the Brownian web has particles starting at all space-time points
$(s,x)$ but we will not need this, and it is enough to establish (\ref{cbmduality}) first for rational
$(x_i)$ and $(a_i)$. The key property is that, almost surely, none of the forward paths cross any of the backwards paths.
(A discrete version of this coupling, that is using simple coalescing simple random walks, is easy
to construct - see the appendix in \cite{stw} - and illustrates this non-crossing property).
From this non-crossing property one sees that the event that $N_t\left((a,b)\right) = 0$ for the forward
coalescing system is almost surely equal to the event that the open interval formed by pair of
backwards particles starting at $a$ and $b$ does not contain any of the initial forwards particles.
The coalescing duality (\ref{cbmduality}) follows immediately, once one notes that
$S_t$ may be replaced by its closure and that annihilating the backwards particles when they meet will
not affect this closure.

The annihilating duality (\ref{abmduality}) follows from (\ref{cbmduality}) and the thinning relation.
Note that thinning a set of $n \geq 1$ elements produces a random subset whose size has a binomial
$B(n,1/2)$ distribution, and also that $E[(-1)^{B(n,1/2)}] =0$.
Then thinning and (\ref{cbmduality}) show that
\begin{eqnarray*}
E^{A}_{\Theta(x_i)} \left[ (-1)^{N_t(I_1 \cup I_3 \cup \ldots \cup I_{2m-1})} \right]
& = & P^{C}_{(x_i)} \left[ N_t(I_1)=N_t(I_3) = \ldots =N_t(I_{2m-1})=0 \right] \\
& = & P^{A}_{(a_i)} \left[ S_t \cap (x_i) = \emptyset \right] \\
& = & E^{A}_{(a_i)} \left[ (-1)^{| S_t \cap \Theta(x_i)|} \right]
\end{eqnarray*}
(where on the right hand side $E^A_{(a_i)}$ is the expectation over the annihilating particle system and over the
independent thinning).
One may then argue by induction on the number $n$ of the initial particles $(x_1, \ldots,x_n)$.
When $n=1$ the above identity reduces to (\ref{abmduality}) for a single particle.
For general $n$ the identity is a mixture of copies of  (\ref{abmduality}) for initial conditions
that are subsets of $(x_i)$. But all but one of the copies will involve
$n-1$ or less particles allowing an inductive proof.
Note also that  (\ref{cbmduality}) also follows from (\ref{abmduality}) - a weighted sum of (\ref{abmduality})
according to the distribution of $\Theta(x_i)$ yields (\ref{cbmduality}).

\noindent
\textbf{Remark. }
Other coalescing duality formulae, such as those in Xiong and Zhou \cite{xiong+zhou}, also
follow from the Brownian web and its dual, but their proof shows that one may also establish them using
the Markov generator duality, as explained in section 4.4 of Ethier and Kurtz, and thus bypass
the Brownian web. In particular this generator technique may be extended to show analogous
dualities for more general spatial motions, where the web construction does not (as yet) exist.
Formally the generator proof shows that the dualities (\ref{cbmduality}) and (\ref{abmduality}) will hold for
instantly reacting continuous Markovian
motions, where the motion on the right hand side must be the image of the motion on the left hand side
under reflection $x \to -x$. Furthermore the maximal entrance laws constructed in the next section
should follow once some moment control is established, which will require some non-degeneracy
of the spatial motion to ensure enough reactions take place.

\subsection{Maximal entrance laws} \label{s2.3}

One may start coalescing systems from infinitely many particles at time zero.
A natural state space for the empirical measure is the set $\mathcal{M}_{LFP}(\R)$
of locally finite point measures on $\R$, which is a closed subset of the space of
locally finite measures under the topology of vague convergence of measures.
The reactions ensure that the point masses only have mass
one, and so we consider the (measurable) subset $\mathcal{M}_0$
of those measures of the form
\[
\mu = \sum_i \delta_{x_i} \quad \mbox{where  $(x_i)$ is locally finite in $\R$ and has disjoint elements.}
\]
(To obtain a process with continuous paths, which does not concern us in this paper, one can quotient
$\mathcal{M}_{LFP}$ by the minimal relation that ensures $\mu + 2 \delta_x \sim \mu + \delta_x$.)

There is a Feller Markov transition kernel $p_t(\mu, d\nu)$ on  $\mathcal{M}_{0}$.
Moreover, there is a maximal entrance law, intuitively starting with one particle at every site (as
in the Brownian web). This can be characterized by passing to the limit in (\ref{cbmduality}) as $(x_i)$
increase to become dense in the real line.  This entrance law, which we denote by
$P^{C}_{\infty}$, has one dimensional distributions satisfying
\begin{equation} \label{cbmentrance}
P^{C}_{\infty} \left[ N_t(I_1)=N_t(I_3) = \ldots =N_t(I_{2m-1})=0 \right]
=
P^{A}_{(a_i)} \left[ \tau < t \right]
\end{equation}
where $\tau$ is the time for complete extinction of the annihilating system.
This characterizes the one dimensional laws on  $\mathcal{M}_{0}$, and these laws are an entrance law
for the Markov transition kernel described above. By the scaling property of Brownian motions we have
$P^{A}_{(Ta_i)}[ \tau < T^2t]$ is independent of $T>0$.
Using (\ref{cbmentrance}) this translates into a scaling for the entrance law
\begin{equation} \label{Scaling}
\mbox{The law of $( T^{-1} X^i_{tT^2}: i \in I_{tT^2})$ is independent of $T>0$ under $P^{C}_{\infty}$.}
\end{equation}
Many suitably spread out and non-degenerate initial conditions are attracted to the maximal entrance law as
$t \to \infty$.  For a large class of initial conditions
$(x_i)$, the law of $(T^{-1} X^i_{T^2t}: i \in I_{T^2t})$
under $P^{C}_{(x_i)}$ converges in distribution, on $\mathcal{M}_{LFP}(\R)$ as $T \to \infty$, to the law of
$(X^i_t:i \in I_t)$ under $P^{C}_{\infty}$. Indeed, using the extension of (\ref{cbmduality})
to countable $(x_i)$, this follows (see the appendix) from
\begin{eqnarray}
&& \hspace{-.8in} P^{C}_{(x_i)} \left[ (T^{-1} X^i_{T^2t}: i \in I_{T^2t})
\cap I_k= \emptyset \; \mbox{for $k=1,3,\ldots,2m-1$} \right]  \nonumber \\
& = & P^{C}_{(x_i)} \left[ N_{tT^2}(TI_1) =  N_{tT^2}(TI_3) = \ldots =  N_{tT^2}(TI_{2m-1}) =0 \right]  \nonumber \\
& = &  P^{A}_{(Ta_i)} \left[S_{tT^2} \cap (x_i) = \emptyset \right]  \nonumber \\
& = &  P^{A}_{(a_i)} \left[ S_t \cap (T^{-1} x_i) = \emptyset \right]  \nonumber \\
& \to & P^{A}_{(a_i)} \left[\tau < t \right]  \nonumber \\
& = & P^{C}_{\infty} \left[ N_t(I_1)=N_t(I_3) = \ldots =N_t(I_{2m-1})=0 \right]. \label{attraction}
\end{eqnarray}
The third equality comes from Brownian scaling and the final equality is
(\ref{cbmentrance}). The convergence holds for deterministic $(x_i)$ for which
$(T^{-1} x_i)$ become dense in any finite interval $[a,b]$ as $T \to \infty$.
A large class of random initial conditions will clearly also work, for example non-zero stationary and spatially ergodic.

For annihilating systems a Markov transition kernel can also be constructed,
using (\ref{abmduality}) and it's extension to countable $(x_i)$ as a means of characterization.
We can define an entrance law
$P^{A}_{\infty}$ for the annihilating system by taking the thinned copy of the entrance law for the
coalescing system. This satisfies the formula
\begin{equation} \label{abmentrance}
E^{A}_{\infty} \left[ (-1)^{N_t(I_1 \cup I_3 \cup \ldots \cup I_{2m-1})} \right]
=
P^{A}_{(a_i)} \left[ \tau < t \right]
\end{equation}
which again determines one dimensional laws on $\mathcal{M}_{0}$ that form an entrance law
for the annihilating system.
The domain of attraction of this entrance law is more delicate. The example in section 3 of
Bramson and Griffeath \cite{bramson+griffeath}
suggests that different approximations to a maximal entrance law may yield different laws at times $t>0$
(their example uses varying intensities of nearby pairs at time zero). For initial conditions that fill the lattice
$\lambda^{-1} \Z$, or that are Poisson
with intensity $\lambda$, the one-dimensional time distributions converge as $\lambda \to \infty$ to those of
the entrance measure, or for a fixed $\lambda$ the large time distribution rescales to those of the entrance
law, by the argument above.

Since we found it difficult to find a full account in the literature, we give, in the
appendix, a brief sketch of the proofs of the results in this subsection.
\section{Proofs} \label{s3}
\subsection{Review of Pfaffians} \label{s3.1}
We give a short summary, targeted at beginners like us,
of the facts we shall use
about Pfaffians (mostly proved in \cite{stembridge} section 2), and of the definition of a Pfaffian point process.
We write $\mbox{Pf}\,(a_{ij}:1 \leq i <j \leq 2n)$ (or just $\mbox{Pf} \, (a_{ij}:i<j)$) for the Pfaffian of the
real anti-symmetric matrix whose elements are $a_{ij}$ for $i<j$.

The determinant of an anti-symmetric matrix of odd order is zero.
Suppose $A$ is an anti-symmetric $2n \times 2n$ matrix.
Then $\det(A)$ is the square of a polynomial of degree $n$ in the matrix elements,
called the Pfaffian of $A$ and written as $\mbox{Pf}\,(A)$.
One can define the Pfaffian as a suitable sum over permutations of products of matrix elements. Indeed,
\begin{equation} \label{permpfaff}
\mbox{Pf}\,(A) = \sum_{\sigma \in \Sigma_{2n}} \mbox{sgn}(\sigma)
a_{i_1,j_1} a_{i_2,j_2} \ldots a_{i_n,j_n}
\end{equation}
where $\Sigma_{2n}$ is the set of permutations $\sigma$ of $\{1,2,\ldots,2n\}$
given by $\sigma(2k-1) = i_k, \; \sigma(2k)=j_k$ for $k=1,\ldots,n$ for which
the choices $(i_k),(j_k)$ satisfy $i_k < j_k$ for all $k$ and $i_1 < i_2 < \ldots < i_n$.
A convenient way to calculate the sign of such a permutation is via {\it crossings}.
The quadruple $i_k,j_k,i_l,j_l$ is called crossed if $i_k < i_l < j_k <j_l$.
Then the sign of $\sigma \in \Sigma_{2n}$
equals $(-1)^M$ where $M$ is the number of crossings.
To visualize these crossings easily one can embed
the integers $1,\ldots,2n$ into the $x$-axis of the plane and join $i_k$ to $j_k$
for each $k$ with a loop in the upper half plane.

It is worth recording the smallest cases:
\[
\mbox{Pf}\, \left(
\begin{array}{cc}
0 & a \\
-a & 0
\end{array}
\right) = a
\qquad
\mbox{Pf}\, \left(
\begin{array}{cccc}
0 & a & b & c \\
-a & 0 & d & e \\
-b & -d & 0 & f \\
-c & -e & -f & 0
\end{array}
\right) = af -be + cd.
\]
The explicit $4 \times 4$ case was used to guess many of the Pfaffian formulae in this paper.

Pfaffians have many similar properties to determinants.
It follows from the definition that
$\mbox{Pf}\,(\lambda_i \lambda_j a_{ij}) =  \mbox{Pf}\,(a_{ij}) \, \prod_k \lambda_k$.
For any $2n \times 2n$ matrix $B$ the product $B^TAB$ is anti-symmetric and
$\mbox{Pf}\,(B^TAB) = \det(B) \mbox{Pf}\,(A)$.
The Pfaffian can be decomposed along a row, or column,
of the matrix. For example if $A$ is a $2n \times 2n$ anti-symmetric matrix it satisfies
the recursion, for any $i \in \{1,2,\ldots,2n\}$,
\begin{equation} \label{recurpfaff}
\mbox{Pf}\,(A) = \sum_{j=1, \, j \neq i}^{2n}  (-1)^{i+j+1}  a_{ij} \mbox{Pf}\,(A^{(i,j)})
\end{equation}
where $A^{(i,j)}$ is the $(2n-2) \times (2n-2)$ submatrix formed by removing the $i$th and $j$th rows and columns.
We will also use a decomposition formula for the Pfaffian of a sum of two
 $2n \times 2n$ anti-symmetric matrices $A$ and $B$, namely
\begin{equation} \label{pfaffsum}
\mbox{Pf}\,(A+B) = \sum_{J} (-1)^{|J|/2} (-1)^{s(J)} \mbox{Pf}\,(A|_{J}) \mbox{Pf}\,(B|_{J^c})
\end{equation}
where: the sum is over all subsets $J \subseteq \{1,2,\ldots,2n\}$ with an even number of terms;
$J^c = \{1,2,\ldots,2n\} \setminus J$; $s(J) =  \sum_{j \in J} j $ (and $s(\emptyset) = 0$); and where
$A|_J$ means the submatrix of $A$ formed by the rows and columns indexed by elements of $J$
(and the Pfaffian of the empty matrix is taken to have value $1$).

Suppose a measurable kernel
\[
K(x,y) = \left( \begin{array}{cc}
K_{11}(x,y) & K_{12}(x,y) \\
K_{21}(x,y) & K_{22}(x,y)
\end{array} \right) \quad \mbox{for $x,y \in \R$}
\]
is anti-symmetric, in the sense $K_{ij}(x,y) = - K_{ji}(y,x)$ for all $i,j \in \{1,2\}$ and $x,y \in \R$.
Suppose it also acts as a kernel for a a bounded operator on $L^2(\R) \oplus L^2(\R)$.
A point process $(X^i: i \in I)$ with $n$-point density functions $\rho^{(n)}(x_1,\ldots,x_n)$ is called
(see Soshnikov \cite{soshnikov})  a {\it Pfaffian point process
with kernel $K$} if $\rho^{(n)}(x_1,\ldots,x_{n})$ is given by the Pfaffian of the
$2n \times 2n$ anti-symmetric matrix formed by the $n^2$ two-by-two matrix entries
$(K(x_i,x_j): i,j =1, \ldots,n)$. The kernel is not uniquely determined.

A very convenient tool for manipulating Pfaffians is the Berezin integral.
We provide arguments that avoid this tool in this paper, and so do not describe the rules for manipulating
these integrals. However they were used repeatedly while exploring these
results, and in the next section we show how the Berezin integral can considerably shorten the argument.
A very readable account of Berezin integrals can be found in  Itzykson and Drouffe \cite{berezin}.
The key property linking the Berezin integral to Pfaffians is (compare with the normalizing determinant for
multi-dimensional Gaussian integrals)
\begin{equation} \label{bint}
\mbox{Pf}\,\left(a_{ij}:i < j\right) = \int d\psi_{2n} \ldots d\psi_{1} e^{- \frac12 \sum_{i,j=1}^{2n} \psi_i a_{ij} \psi_j}.
\end{equation}
\subsection{Proof of Theorem \ref{T3}, the product moment Pfaffians} \label{s3.2}
We start with the product moment, defined for bounded measurable $g:\R \to \R$ and
disjoint $(x_i)$ by
\[
m^{(n)}_t(x_1,\ldots,x_{n}) =  E^{A}_{(x_1,\ldots,x_{n})} \left[ \prod_{i \in I_t} g(X^i_t) \right]
\]
where the product over an empty set, occurring when all the particles have been annihilated,
is defined to have value $1$. Note that $m^{(1)}_t(x)$ is given by the Brownian semigroup applied to $g$.
We also set $m^{(0)}\equiv 1$.
To show that $m^{(2n)}$ is given by a Pfaffian, we shall give a p.d.e. derivation similar in spirit to
that showing the Karlin McGregor formula for the transition
density for non-intersecting Brownian motions is given by a determinant.
%However, one can in fact prove that the various Pfaffians in this section hold true in different
%ways and starting from different places.
%Here, the Pfaffian is of a matrix where variables $x=(x_1,\ldots,x_n)$ are the initial conditions. Later
%duality fromulae will give more useful Pfaffians of matrices where the variables will be in the the
%functional applied at time $t$.

Let $V_n \subseteq \R^{n}$ be the open cell $\{x:x_1 < x_2 < \ldots < x_{n}\}$.
On $(0,\infty) \times V_n$ the function $m^{(n)}_t(x)$ solves the heat equation, and we
must examine the boundary conditions.
For $n \geq 2$ and when $g$ is bounded and continuous, the functions $m^{(n)}$
are continuous on $[0,\infty) \times V_n$ and extend to a continuous function in
$C((0,\infty) \times \overline{V}_n)$. There are lots of pieces to the boundary of $V_n$, but the most important
are the faces $F_{i,n}$ defined by $x_i = x_{i+1}$ and where the remaining $x_k$ are disjoint. On $F_{i,n}$ the continuous
extension agrees with
the lower order moment $m^{(n-2)}(x^{(i,i+1)})$, where $x^{(i,j)} \in \R^{n-2}$ is the $(n-2)$-tuple
formed by removing $x_i$ and $x_j$ from $(x_1,\ldots,x_n)$. This can be seen by showing that near the boundary
the hitting time between particles starting at $x_i$ and $x_{i+1}$ is likely to occur before any other
collision and before time $t$. On other parts of the boundary the extension agrees with other lower moments.

The system of heat equations for $(m^{(n)}: n=1,2,\ldots)$
\[
\left\{ \begin{array}{lcll}
\frac{\partial}{\partial t} m^{(n)}_t(x) & = & \Delta m^{(n)}_t(x) & \mbox{on $(0,\infty) \times V_n$,} \\
m^{(n)}_t(x) & = & m^{(n-2)}_t(x^{(i,i+1)}) & \mbox{for $x \in F_{i,n}$ and $i=1,\ldots,n-1$,} \\
m^{(n)}_0(x) & = & \prod_{i=1}^n g(x_i) & \mbox{for $x \in V_n$,}
\end{array} \right.
\]
forms a closed system, in that each equation has boundary
conditions formed by equations of lower order. Note that, typically, the initial condition does not match
the boundary conditions. Taking $g$ bounded and smooth, the system has unique solutions in
$C^{1,2}([0,\infty) \times V_n) \cap C((0,\infty) \times \overline{V}_n)$. It is enough to specify
boundary conditions only on
each face $F_{i,n}$ - the Feynman-Kac formula makes it clear that the other parts of the boundary
of $V_{n}$ do not affect the value of $m^{(n)}$.

To establish the Pfaffian  (\ref{evenpf}) stated in Theorem \ref{T3}, it is enough, by an approximation argument, to treat the case where $g$
is smooth.
We shall prove (\ref{evenpf}) by showing the Pfaffian $ \mbox{Pf} \, \left(m^{(2)}_t(x_i,x_j): 1 \leq i < j \leq 2n \right)$
solves the system of heat equations above.
Note that (\ref{evenpf}) holds when $t=0$ since
\[
m^{(2n)}_0(x_1,\ldots,x_{2n}) = \prod_{i=1}^{2n} g(x_i) = \mbox{Pf}\, \left( g(x_i) g(x_j)\,:\, i < j \right).
\]
The Pfaffian is a finite sum of product terms (see (\ref{permpfaff})) of the form
\[
\mbox{sgn}(\sigma) m^{(2)}_t(x_{i_1},x_{j_1}) m^{(2)}_t(x_{i_2},x_{j_2}) \ldots m^{(2)}_t(x_{i_n},x_{j_n})
\]
where $\sigma$ is a permutation given by $\sigma(2k-1) = i_k, \; \sigma(2k)=j_k$ for $k=1,\ldots,n$.
Since $m^{(2)}_t(x,y)$ satisfies the heat equation on $[0,\infty) \times \{x < y\}$, each product
term lies in $C^{1,2}([0,\infty) \times V_{2n})$ and
satisfies the heat equation on $[0,\infty) \times V_{2n}$. Since $m^{(2)}_t(x,y)$ extends continuously to
$(0,\infty) \times \{(x,y): x \leq y\}$, the Pfaffian extends continuously to $(0,\infty) \times \overline{V}_{2n}$.
By the uniqueness for the system of heat equations, it remains to check that the Pfaffian satisfies the required
boundary conditions on each face $F_{i,2n}$ which will complete the proof of Theorem \ref{T3}.

We show the argument for the face $F_{1,2n}$ where $x_1=x_2$ (other faces are similar).
We may argue inductively, and suppose
that $m^{(k)}$ is given by the Pfaffian for $k=0,2,\ldots,2n-2$.  Our quickest proof is
using the representation (\ref{bint}) in terms of Berezin integrals. This gives
\begin{eqnarray*}
&& \hspace{-.3in}  \left.  \mbox{Pf}\, \left( m^{(2)}_t(x_i,x_j)\,: \, i <j \right) \right|_{x_1=x_2} \\
& = &  \left. \int d\psi_{2n} \ldots d\psi_{1} e^{- \frac12 \sum_{i,j=1}^{2n} \psi_i m^{(2)}_t(x_i,x_j) \psi_j} \right|_{x_1=x_2} \\
& = &   \int d\psi_{2n} \ldots d\psi_{1} e^{- \frac12 \sum_{i,j=3}^{2n} \psi_i m^{(2)}_t(x_i,x_j) \psi_j}
 e^{- (\psi_1 + \psi_2) \sum_{k=3}^{2n} m^{(2)}_t(x_1,x_k) \psi_k }.
\end{eqnarray*}
The sum $M = \sum_{k=3}^{2n} m^{(2)}_t(x_1,x_k) \psi_k$ is independent of $\psi_1$ and $\psi_2$ and the
 $d\psi_{2} d\psi_{1}$ integral becomes (using the rules for Berezin integrals)
\[
\int d\psi_{2} d\psi_{1} e^{- (\psi_1 + \psi_2) M} = \int d\psi_{2} d\psi_{1} (1-\psi_2 \psi_1)(1-(\psi_1+\psi_2)M) = 1.
\]
This simplification leaves only $ \int d\psi_{2n} \ldots d\psi_{3}
e^{- \frac12 \sum_{i,j=3}^{2n} \psi_i m^{(2)}_t(x_i,x_j) \psi_j}$
which is the Berezin integral for $m^{(2n-2)}(x_3,\ldots,x_{2n})$.

An argument that avoids Berezin integrals is as follows.
Using the recursive relation for Pfaffians (\ref{recurpfaff}) we see that the Pfaffian in (\ref{evenpf})
equals
\[
\sum_{k=2}^{2n} (-1)^{k} m^{(2)}_t(x_1,x_k) m^{(2n-2)}_t(x^{(1,k)}).
\]
Since $m^{(2)}_t(x_1,x_2)$ extends to the function $1$ on $x_1=x_2$, it remains only to check that
\begin{equation} \label{temp10}
\sum_{k=3}^{2n} (-1)^{k} m^{(2)}_t(x_1,x_k) m^{(2n-2)}_t(x^{(1,k)})
\end{equation}
vanishes when $x_1=x_2$ and $t>0$. But this follows from expressing $m^{(2n-2)}$ using
(\ref{permpfaff}). Indeed, fix $j,k \geq 3$. Then for an expression of the form
\[
m^{(2)}_t(x_1,x_k) m^{(2)}_t(x_2,x_j) m^{(2)}_t(x_{i_2},x_{j_2}) \ldots m^{(2)}_t(x_{i_{n-1}},x_{j_{n-1}})
\]
arising from the $k$th term in (\ref{temp10}), where
$\{i_2,j_2,\ldots,i_{n-1},j_{n-1}\} = \{3,4,\ldots,2n\} \setminus \{j,k\}$,
 there is a corresponding term
\[
m^{(2)}_t(x_1,x_j) m^{(2)}_t(x_2,x_k) m^{(2)}_t(x_{i_2},x_{j_2}) \ldots m^{(2)}_t(x_{i_{n-1}},x_{j_{n-1}})
\]
arising from the $j$th term in (\ref{temp10}). These terms agree on $x_1=x_2$ and
a careful check of the signs of the permutations, and the factors $(-1)^j$ and $(-1)^k$ in
(\ref{temp10}), shows they will cancel.
One way to do this check is to compare the sign of
\[
\sigma = \left( \begin{array}{cccccccccc} 2 & 3 & 4 & 5 & \ldots & k-1 & k+1& \ldots & 2n-1 & 2n \\
                                                 2 & j & i_2 & i_3 &\ldots & \ldots & \ldots &\ldots & i_{n-1}& j_{n-1}
\end{array} \right)
\]
with that of
\[
\sigma' = \left( \begin{array}{cccccccccc} 2 & 3 & 4 & 5 & \ldots & j-1 & j+1& \ldots & 2n-1 & 2n \\
                                                 2 & k & i_2 & i_3 & \ldots & \ldots & \ldots &\ldots & i_{n-1}& j_{n-1}
\end{array} \right)
\]
by counting crossings. The loop joining $2$ to $j$ in $\sigma$ must be replaced by a loop joining $2$ to $k$ in
$\sigma'$. This may affect crossings with any of the loops emanating from sites between $j$ and $k$,
and will do so unless a pair of them are joined to each other.
There are $|k-j|-1$ sites between $j$ and $k$ so it will change the parity of the number of crossings exactly
when $|k-j|$ is even.

\noindent
\textbf{Remark 1.}
For odd moments there is also a Pfaffian representation, namely,
when $x_1 < x_2 < \ldots < x_{2n-1}$,
\begin{equation} \label{oddpf}
m^{(2n-1)}_t(x_1,\ldots,x_{2n-1}) = \mbox{Pf}\, \left( m^{(2)}_t(x_{i},x_{j}): 0 \leq i < j \leq 2n-1 \right)
\end{equation}
where we adopt the convention that $m^{(2)}_t(x_0,x_k) = m^{(1)}_t(x_k)$.
This Pfaffian involves a linear combination of terms of the form
\[
\mbox{sgn}(\sigma) m^{(1)}_t(x_{j_1}) m^{(2)}(x_{i_2},x_{j_2}) \ldots m^{(2)}(x_{i_n},x_{j_n})
\]
which again shows that it solves the heat equation when $[0,\infty) \times V_{2n-1}$.
The recursive Pfaffian relation gives
\[
m^{(2n-1)}_t(x) = \sum^{2n-1}_{k=1} (-1)^k m^{(1)}_t(x_k) m^{(2n-2)}_t(x^{(k)}).
\]
Expanding the Pfaffian along its first row using (\ref{recurpfaff}) we obtain
for $x = (x_1,\ldots,x_{2n-1}) \in V_{2n-1}$
\begin{equation} \label{temp20}
\mbox{Pf}\, \left( m^{(2)}_t(x_{i},x_{j}): 0 \leq i < j \leq 2n-1 \right)
= \sum_{k=1}^{2n-1} (-1)^{k+1} m^{(1)}_t(x_k) m^{(2n-2)}_t(x^{(k)})
\end{equation}
where we again write superscripts $x^{(i,j,\ldots)}$ to mean that we remove the indicated co-ordinates.
The terms with $k=1$ and $k=2$ cancel on the face $F_{2n-1,1}$ where $x_1=x_2$.
Moreover on this face, for $k \geq 3$, $m^{(2n-2)}_t(x^{(k)}) = m^{(2n-4)}_t(x^{(1,2,k)})$
so that the Pfaffian in (\ref{temp20}) becomes
\[
\sum_{k=3}^{2n-1} (-1)^{k+1} m^{(1)}_t(x_k) m^{(2n-4)}_t(x^{(1,2,k)}).
\]
But this is the decomposition of $m^{(2n-3)}_t(x^{(1,2)})$ along the first row, and
this shows the boundary conditions are correct on $F_{2n-1,1}$. Other faces are similar.

\noindent
\textbf{Remark 2.}  Since our proof relies only on uniqueness
for the underlying system of heat equations, the extension of these product moment Pfaffians
to more general spatial motions looks quite straightforward, for example to Markovian spatial motions
that are continuous and suitably non-degenerate. The Pfaffians in the next section
would then also follow for these more general motions, just by algebraic manipulation, once maximal
entrance laws characterized by (\ref{cbmentrance}) and (\ref{abmentrance}) are established.
%We anticipate that a detailed
%treatment will be given in the Warwick thesis \cite{yip}.
%
%
%
\subsection{Proof of Theorem \ref{T2}, the Pfaffian point process kernel} \label{s3.3}
Fixing $a_1< \ldots <a_{2m}$ and choosing $ g(x) = (-1)^{\sum_i \chi(x \leq a_i)}$ in (\ref{evenpf})
we see that both sides of the duality (\ref{abmduality}) are Pfaffians in the variables $(x_i)$.
Choosing $g=0$, recalling that an empty product takes the value $1$, we see that $P^{A}_{(x_i)}[\tau <t]$ is a Pfaffian.
The entrance law dualities (\ref{cbmentrance}) and (\ref{abmentrance}) show that
\[
P^{C}_{\infty} \left[ N_t(I_1)=N_t(I_3) = \ldots =N_t(I_{2m-1})=0 \right]
= E^{A}_{\infty} \left[ (-1)^{N_t(I_1 \cup I_3 \cup \ldots \cup I_{2m-1})} \right]
\]
are Pfaffians in the variables $(a_i)$. The entries in this last Pfaffian are explicit since
\[
 P^{C}_{\infty} \left[ N_t \left((a_j,a_k)\right) =0 \right]
= E^{A}_{\infty} \left[ (-1)^{N_t \left((a_j,a_k)\right)} \right]
=  P^{A}_{(a_i,a_j)}\left[\tau <t \right]
\]
are all equal to (by Brownian hitting probabilities)
\begin{equation}
\mbox{$ F\left( t^{-1/2}(a_j-a_i)\right) $
where $F(x) = \pi^{-1/2} \int^{\infty}_x \exp(-y^2/4) \, dy$.} \label{Amatrix}
\end{equation}
We switch dummy variables for the rest of this section, taking $x_1<x_2< \ldots <x_{2n}$
and $I_k=(x_k,x_{k+1})$ so that we start from
\begin{equation} \label{cbmpf}
P^{C}_{\infty} \left[ \mbox{$N_t(I_k)=0$ for $k=1,3,\ldots,2n-1$} \right]
= \mbox{Pf} \left(F\left(t^{-1/2}(x_j-x_i)\right)\,:\, i < j \right).
\end{equation}

To prove Theorem \ref{T2}, we shall identify the Pfaffian point process kernel
by differentiating the empty interval Pfaffian (\ref{cbmpf}) above.  By scaling we may take
$t=1$. Differentiate the identity (\ref{cbmpf}) for $t=1$ in the variables $x_1,x_3,\ldots,x_{2n-1}$. The left hand side
becomes, formally,
\[
E^{C}_{\infty} \left[ N_1(dx_1) \, N_1(dx_3) \ldots N_1(dx_{2n-1})
\I(\mbox{$N_1(I_k)=0$ for $k=1,3,\ldots,2n-1$}) \right].
\]
Letting $x_{2l} \downarrow x_{2l-1}$ for $l=1,\ldots,n$ we reach the
$n$-point density $\rho^{(n)}_1(x_1,x_3,\ldots,x_{2n-1})$.
In the appendix \ref{s4.3} we give more details verifying the formal
differentiation above is valid, by using distributional derivatives.

On the right hand side of (\ref{cbmpf}) we will also
differentiate in the variables $x_1,x_3,\ldots,x_{2n-1}$.
Note that each product term in the permutation expansion (\ref{permpfaff}) of the Pfaffian contains
exactly one element that involves the variable $x_1$. So differentiating in $x_1$ leads to
a similar permutation expansion, but where all the terms that involve $x_1$ have been differentiated.
Repeating this argument, differentiating in $x_1,x_3,\ldots,x_{2n-1}$ yields the Pfaffian where
each term in the matrix has been differentiated in the variables $x_1,x_3,\ldots,x_{2n-1}$
that is where the $2 \times 2$ block formed by the rows
$2j-1,\,2j$ and columns $2k-1,\,2k$ is given by
\[
\left( \begin{array}{cc}
- F''(x_{2k-1}-x_{2j-1}) & -F'(x_{2k}-x_{2j-1})  \\
F'(x_{2k-1}-x_{2j}) & \sgn(x_{2k}-x_{2j}) \,   F(x_{2k}-x_{2j})
\end{array} \right)
\]
when $j \leq k$. (Note that $F''$ is an odd function and so no $\sgn$ is needed in the $2j-1, 2k-1$ entry.)
 Letting $x_{2l} \downarrow x_{2l-1}$ for $l=1,\ldots,n$ we obtain
the kernel $K$ stated in Theorem \ref{T2}. The decay in $F,F',F''$ implies that $K$ acts as a
suitable bounded operator. The scaling relation (\ref{Scaling}) implies that the kernel of the distribution
 time $t$ is $t^{-1/2} K(x t^{-1/2},yt^{-1/2})$.

\noindent
\textbf{Remark 1.}
An alternative starting point, used by ben Avraham et al. (see \cite{avraham1}, \cite{avraham2}),
is to show the empty interval probabilities
$P^{C}_{\infty} [ N_t(I_1)=N_t(I_3) = \ldots =N_t(I_{2m-1})=0]$
satisfy heat equations in the variables $(x_i)$, though
the connection with Pfaffians does not seem to have been noted.

\noindent
\textbf{Remark 2.} The linear ordering of particles seems to be crucial. For Brownian particles
on a one-dimensional torus, there is an extra boundary condition where $x_{2n}$ may hit $x_1$ by going 'the other way'
around the torus, and this is not satisfied by the Pfaffians.

\noindent
\textbf{Remark 3.} Differentiating in the variables $x_2, \ldots,x_{2n}$ instead leads to the
alternative dual kernel
\[
K(x,y) = \left( \begin{array}{cc}
\sgn(y-x) F(|y-x|) & F'(y-x) \\
-F'(y-x) &  -F''(y-x)
\end{array} \right)
\]

\noindent
\textbf{Remark 4.}
Starting from the Pfaffian (\ref{cbmpf}) certain other probabilities can, by
algebraic manipulation, also be expressed as Pfaffians.
We give three examples, leaving details of the derivations to the appendix. In
each case $F$ is the $2n \times 2n$ anti-symmetric matrix with elements
$F_{ij} = P^{C}_{\infty} [N_t \left((x_j,x_k)\right)=0] = F(t^{-1/2}(x_j-x_i))$ as
in (\ref{Amatrix}).
\begin{itemize}
\item
Let $I=I_{2n}$ be the $2n \times 2n$ anti-symmetric matrix with entries $1$ above the diagonal. Then
\begin{equation} \label{cbmpf2}
P^{C}_{\infty} \left [\mbox{$N_t(I_k) > 0 $ for $k=1,\ldots,2n-1$} \right] =
\mbox{Pf}\,(I-F).
\end{equation}
The annihilating analogue of this is
\begin{equation} \label{cbmpf2.5}
P^{A}_{\infty} \left [\mbox{$N_t(I_k )$ is odd for $k=1,2,\ldots,2n-1$} \right] =
2^{1-2n} \mbox{Pf}\,(I-F).
\end{equation}
\item Let $O=O_{2n}$ be the $2n \times 2n$ anti-symmetric matrix formed by $n$ copies of the
$2 \times 2$ matrix {\small $\left( \begin{array}{ll} 0 & 1 \\ -1 & 0 \end{array} \right)$}
down the diagonal and zeros elsewhere. Then
\begin{equation} \label{cbmpf3}
P^{C}_{\infty} \left [\mbox{$N_t \left(I_k)\right) > 0 $ for $k=1,3,5,\ldots,2n-1$} \right]
= \mbox{Pf}\,(O-F).
\end{equation}
Again there is an annihilating analogue.
\item Let $\hat{O}=\hat{O}_{2n}$ be the $2n \times 2n$ anti-symmetric matrix with entries
\[
\hat{O}_{ij} = \left\{ \begin{array}{cl}
+1 & \mbox{if $i=2,4,\ldots,2n-2$ and $j=i+1$,} \\
-1 & \mbox{if $j=2,4,\ldots,2n-2$ and $i=j+1$,} \\
0 & \mbox{otherwise.} \end{array} \right.
\]
Note that $\hat{O}$ also has copies of the  $2 \times 2$ matrix
{\small $\left( \begin{array}{ll} 0 & 1 \\ -1 & 0 \end{array} \right)$}
in some places down the diagonal and zeros elsewhere. Then
\begin{eqnarray}
&& \hspace{-.3in}
P^C_{\infty} \left[ \mbox{$N_t(I_k)=0$ for $k=1,3,\ldots,2k-1$ and  $N_t(I_k)>0$ for
$k=2,4,\ldots,2k-2$} \right] \nonumber \\
& = &  \mbox{Pf}\,(F-\hat{O}). \label{cbmpf4}
\end{eqnarray}
\end{itemize}

\noindent
\textbf{Remark 5.}
Suppose $(M_{xy}:x \leq y)$ is a bounded continuous field that satisfies
$E[M_{x_1x_2} M_{x_3x_4}] = \mbox{Pf}\,(E[M_{x_ix_j}]: 1 \leq i < j \leq 4)$
for $x_1< \ldots <x_4$. Then by continuity $E[M_{xx}M_{xx}] = \mbox{Pf}\,(E[M_{xx}])
= E[M_{xx}]^2$ and so $M_{xx}$ must be deterministic. This imposes a restriction
on the class of correlation functions admitting a Pfaffian representation.
\subsection{Proof of Theorem \ref{T1}, the asymptotics for $\rho^{(2n)}_t$} \label{s3.4}
%
%In this section we give the proof of Theorem \ref{T1}.
We work throughout under the entrance measure $P^{C}_{\infty}$.
By thinning the corresponding density for
annihilating systems differs only by a multiplicative factor $2^{-n}$.
The $n$-point density function $\rho^{(n)}_t(x)$ is a Lebesgue density for the measure
$E^{C}_{\infty}[N_t(dx_1) \ldots N_t(dx_n)]$ on $V_{n}$.
The existence of this density, defined almost everywhere, and the simple bound
\begin{equation} \label{simplebound}
\rho_t^{(n)} (x) \leq C_n t^{-n/2} \quad \mbox{for all $t>0$ and $x \in V_{n}$}
\end{equation}
is discussed in \cite{mrtz}. Furthermore there we established the following upper bound:
for all $L>0$ there exists $C_L < \infty$ so that
\begin{equation} \label{oldasympt}
\rho_t^{(n)}(x_1,\ldots,x_n) \leq C_L t^{-\frac{n}{2}-\frac{n(n-1)}{4}} \prod_{1 \leq i < j \leq n} |x_i-x_j|
\quad \mbox{for all $t>0$ and $|x_i| \leq Lt^{1/2}$}.
\end{equation}

As $t \to \infty$ the entries in the Pfaffian for $\rho^{(n)}_t$ are of the form $F,\,F',\,F''$ evaluated at
points $t^{-1/2}(x_j-x_i)$ close to zero. One may approximate these by using the Taylor expansion for
$F(z)$ at small values of $z$. However, considerable cancellation occurs in the many terms of the Pfaffian and it is not
immediately clear how to read off the leading asymptotic decay in $t$. Indeed the following argument shows
at $F$ needs to be expanded to a large number of terms to obtain the correct answer.

We shall analyze first a modified density function $\tilde{\rho}^{(2n)}_t(x)$ for
$x \in V_{2n}$, which is a density for the measure
\[
E^{C}_{\infty} \left[ N_t(dx_1) \ldots N_t(dx_{2n}) \chi(\mbox{$N_t(I_k)=0$ for $k=1,3,\ldots,2n-1$}) \right] \\
\]
(where we recall that $I_k=(x_k,x_{k+1})$).
%If we label the particles alive at time $t$ as $(X^i_t:i \in \Z)$ with the convention that
%$X^i_t < X^{i+1}_t$ for all $i$, then $\tilde{\rho}^{(2n)}_t(x)$ satisfies, for measurable $f \geq 0$,
%\begin{eqnarray*}
%&& \hspace{-.4in} \int_{V_{2n}} f(x_1,\ldots,x_{2n}) \tilde{\rho}^{(2n)}_t(x_1,\ldots,x_{2n}) dx_1 \ldots dx_{2n} \\
%& = &  E^{C}_{\infty} \left[ \int_{V_{2n}}  f(x_1,\ldots,x_{2n})
 %\chi(\mbox{$N_t(I_k)=0$ for $k=1,3,\ldots,2n-1$})
 %N_t(dx_1) \ldots N_t(dx_{2n})  \right] \\
%& = & E^{C}_{\infty} \left[ \sum_{i_1} \sum_{i_2>i_1+1} \ldots \sum_{i_n>i_{n-1}+1}
%f(X^{i_1}_t,X^{i_1+1}_t,X^{i_2}_t,X^{i_2+1}_t, \ldots,X^{i_n}_t,X^{i_n+1}_t) \right].
%\end{eqnarray*}
%The usual $n$-point density $\rho^{(n)}_t(x)$ satisfies
%\begin{eqnarray*}
%  && \hspace{-.4in} \int_{V_{2n}} f(x_1,\ldots,x_{2n}) \rho^{(n)}_t(x_1,\ldots,x_{2n}) dx_1 \ldots dx_{n} \\
%&  = & E^{C}_{\infty} \left[ \sum_{i_1} \sum_{i_2>1_i} \ldots \sum_{i_{n}>i_{n-1}}
%f(X^{i_1}_t,X^{i_2}_t,\ldots, X^{i_{n}}_t) \right].
%\end{eqnarray*}
We claim that
\begin{equation} \label{tildeclaim}
\tilde{\rho}_t^{(2n)}(x_1,\ldots,x_{2n}) = (4 \pi t^2)^{-n/2} \,
\mbox{Pf}\, \left (\phi\left( (x_j-x_i)/t^{1/2} \right): 1 \leq i < j \leq 2n\right)
\end{equation}
where $ \phi(z) = z \exp (-z^2/4)$. This follows formally, as in section \ref{s3.3}, by differentiating
(\ref{cbmpf}) in all variables $x_1,x_2,\ldots,x_{2n}$, and using that,
\[
  \partial_{x_1} \ldots \partial_{x_{2n}} \mbox{Pf}\, \left(F(x_j-x_i)\,:\,i<j \right)
= \mbox{Pf}\,\left(-(4\pi)^{-1/2} \phi(x_j-x_i)\right)
\]
(which follows from differentiating each term in the permutation expansion (\ref{permpfaff}) of the Pfaffian).
We give more details in the appendix \ref{s4.3}.

The advantage of the representation (\ref{tildeclaim}) is that
it is a Pfaffian all of whose entries are of the form $f(x_i-x_j)$ for a single function $f$, and
this allows us to apply the following lemma, proved at the end of this section,
that gives an expansion for a Pfaffian whose entries are close to the zero of an odd function.
\begin{Lemma} \label{olegslemma}
Let $\phi:\R \to \R$ be an odd function that is analytic at zero. Then for any $n \geq 1$
there exist $\epsilon(n,\phi)>0$ and $C(n,\phi) < \infty$ so that for
$y \in V_n$ with $|y| \leq \epsilon(n,\phi)$
\[
\mbox{Pf}\, \left(\phi(y_j - y_i): 1 \leq i < j \leq 2n \right) =  \mbox{Pf}\,(J^{(2n)}(\phi) + R^{(2n)}(y)) \,
\prod_{1 \leq i<j \leq 2n} (y_j-y_i)
\]
where $J^{(2n)}(\phi) $ is the constant anti-symmetric matrix with entries
\begin{equation}
J^{(2n)}_{ij}(\phi) = (-1)^{j-1} \frac{1}{(i-1)! (j-1)!} \frac{d^{i+j-2} \phi}{d x^{i+j-2}} (0) \quad \mbox{for $1 \leq i <j \leq 2n$,}
\label{J}
\end{equation}
and the remainder $R^{(2n)}(y)$ is a anti-symmetric matrix satisfying
\[
\left| R^{(2n)}_{ij}(y) \right| \leq C(n,\phi) |y| \quad \mbox{for all $i,j$ and $|y| \leq \epsilon(n,\phi)$.}
\]
\end{Lemma}
We apply this lemma to the Pfaffian in (\ref{tildeclaim}) with $\phi(z) = z e^{-z^2/4}$ and with
$y= t^{-1/2} x$ for $t$ large enough.
Expanding the Pfaffian $\mbox{Pf}\,(J^{(2n)}(\phi) + R^{(2n)}(t^{-1/2}x))$ using
(\ref{pfaffsum}) we find only one term, namely $\mbox{Pf}\,(J^{(2n)}(\phi))$, that does not
decay as $t \to \infty$. This shows that
\[
\lim_{t \to \infty} t^{n^2+\frac{n}{2}} \tilde{\rho}^{(2n)}_t(x_1,\ldots,x_{2n})
= (4\pi)^{-n/2} \, \mbox{Pf}\,\left(J^{(2n)}(\phi)\right) \, \prod_{1 \leq i<j \leq 2n}
(x_j-x_i).
\]
To obtain the same estimate for $\rho^{(2n)}$ we estimate the difference as follows.
\begin{eqnarray*}
0 & \leq & \rho^{(2n)}_t(x_1,\ldots,x_{2n}) - \tilde{\rho}^{(2n)}_t(x_1,\ldots,x_{2n}) \\
%& = & E^{C}_{\infty} \left[ N_t(dx_1) \ldots N_t(dx_{2n})
%\left( 1-\chi(\mbox{$N_t(I_k)=0$ for $k=1,3,\ldots,2n-1$}) \right) \right] \\
& = & E^{C}_{\infty} \left[ N_t(dx_1) \ldots N_t(dx_{2n})
\chi(\mbox{$N_t(I_k)>0$ for some $k=1,3,\ldots,2n-1$}) \right] \\
& \leq & \sum_{k=1}^{2n-1} E^{C}_{\infty} \left[ N_t(dx_1) \ldots N_t(dx_{2n}) N_t(I_k) \right] \\
& = & \sum_{k=1}^{2n-1} \int_{I_k} \rho^{(2n+1)}_t(x_1,\ldots,x_k,z,x_{k+1},\ldots,x_{2n}) dz.
\end{eqnarray*}
Each term in this sum is of a smaller order in $t$ by (\ref{oldasympt}).

Examination of the proof shows that we need not let the values of $x_1,\ldots,x_{2n}$ be fixed, and that
in fact we may take the supremum over any positions
$(x_{i}(t))$ provided that $\sup_i |x_i(t)|t^{-1/2} \to 0$ as $t \to \infty$.

\noindent
\textbf{Proof of Lemma \ref{olegslemma}.}
Let $\Phi$ be the $2n \times 2n$ anti-symmetric matrix with entries given by $\Phi_{ij} = \phi(y_j-y_i)$. The aim is to show,
for small $y$, that
\[
\Phi = V^T(J+R)V
\]
 where $J$ and $R$ are as in the lemma (with $n$ fixed and suppressed)
and $V$ is the $2n \times 2n$ Vandermond matrix
given by $V_{ij} = y_j^{i-1}$. Since $\det(V) = \prod_{1 \leq i <j \leq 2n} (y_j-y_i)$,
the conclusion then holds from $\mbox{Pf}\, ( V^T(J+R)V) = \det(V) \mbox{Pf} \, (J+R)$

For small $|y|$ we expand by analyticity (writing $\phi^k(0)$ for the $k$th derivative of $\phi$ at zero)
\begin{eqnarray}
\Phi_{ij}
& = & \sum_{n=0}^{\infty}\frac{1}{n!} \phi^n(0) (y_j-y_i)^n \nonumber \\
& = & \sum_{n=0}^{\infty}  \sum_{k=0}^n \frac{1}{k! \, (n-k)!} \phi^n(0) y^k_j (-y_i)^{n-k} \nonumber \\
& = & \sum_{k,l=0}^{\infty}  \frac{1}{k! \, l!} \phi^{k+l}(0) y^k_j (-y_i)^{l} \nonumber \\
& = & \sum_{k,l=1}^{\infty} y_i^{l-1} y_j^{k-1} J_{lk} \label{newton}
\end{eqnarray}
where we have rearranged using $l=n-k$ in the penultimate equality. Note that
\[
(V^T J V)_{ij} = \sum_{k,l=1}^{2n} V_{li} V_{kj} J_{lk} = \sum_{k,l=1}^{2n} y_i^{l-1} y_j^{k-1} J_{lk}.
\]
It remains to re-express the remaining terms in (\ref{newton}) as the desired remainder.

Recall the symmetric polynomials $\sigma_k^{2n}(y)$ defined for $y \in \R^{2n}$ by
\begin{equation} \label{symmetric}
\prod_{k=1}^{2n} (y_k - \lambda) = \sum_{k=0}^{2n} (-1)^k \sigma_k^{2n}(y) \lambda^{2n-k}.
\end{equation}
Note that $\sigma^{2n}_k$ is a polynomial of order $k$.
Since $\sigma^{2n}_0 \equiv 1$ we may choose $\lambda = y_i$ to see that
\[
0 = y_i^{2n} + \sum_{k=1}^{2n}  (-1)^k \sigma_k^{2n}(y) y_i^{2n-k} \qquad \mbox{for $i=1,\ldots,2n$.}
\]
Multiplying by $y_i^p$ we see that
\begin{equation} \label{pretau}
y_i^{p+2n} = \sum_{k=1}^{2n}  (-1)^{k+1} \sigma_k^{2n}(y) y_i^{p+2n-k} \qquad \mbox{for $i=1,\ldots,2n$ and $p=0,1,\ldots$.}
\end{equation}
By iterating this we may express $y_i^{p+2n}$ for $p \geq 0$ as a mixture of $1,y_i,y_i^2, \ldots,y_i^{2n-1}$, as follows:
\begin{equation} \label{tau}
y_i^{p+2n} = \sum_{k=1}^{2n} \tau^{2n,p+2n}_k(y) y_i^{k-1}  \qquad \mbox{for $i=1,\ldots,2n$ and $p=0,1,\ldots$}
\end{equation}
where $\tau^{2n,p+2n}_k(y)$ is a polynomial of order $p+2n-k+1$.
Using this substitution in the remaining terms of (\ref{newton}), that is where $k$ or $l$ is at least $2n+1$, we find
(formally) that
\[
 \left( \sum_{k,l=2n+1}^{\infty} +  \sum_{k=1}^{2n} \sum_{l=2n+1}^{\infty} +  \sum_{l=1}^{2n} \sum_{k=2n+1}^{\infty} \right)
y_i^{l-1} y_j^{k-1} J_{lk} = \sum_{p,q=1}^{2n} y_i^{p-1} y_j^{q-1} R_{pq}(y)
\]
where
\begin{eqnarray}
R_{pq} (y) & = &  \sum_{k,l=2n+1}^{\infty} \tau^{2n,l-1}_p(y) J_{kl} \tau^{2n,k-1}_q(y)  \nonumber \\
&& +  \sum_{k=1}^{2n} \sum_{l=2n+1}^{\infty}  \tau^{2n,l-1}_p(y) J_{ql}
+  \sum_{l=1}^{2n} \sum_{k=1}^{\infty}   \tau^{2n,k-1}_q(y) J_{kp}. \label{Reqn}
\end{eqnarray}
Note the lowest order of the polynomial entries in the terms for $R_{pq}$ is of order $1$. In the appendix \ref{s4.4} we
check that this rearrangement of (\ref{newton}) is valid when $|y|$ is suitably small and that the required
error bound $|R_{pq}(y)| \leq C(n,\phi) |y|$ holds.

\noindent
\textbf{Thanks.} We would like to thank our colleague Dmitriy Rumynin for advice on the use of symmetric polynomials.
\section{Appendix} \label{s4}
\subsection{Details for section \ref{s2.3}}  \label{s4.1}
We give a few details on (one approach to) the results surveyed  in section \ref{s2.3}.
For coalescing systems one can use monotonicity, adding initial particles one by one, to construct the
infinite system. This is not available for annihilating systems, so we sketch a weak convergence argument
that applies to both systems.

One can control moments by bounds on the $n$-point density function. Indeed $\rho^{(n)}_t(x)$, the density
for the measure $E^{C}_{(x_i)}[N_t(dx_1) \ldots N_t(dx_n)]$ on $V_n$, depends on the initial condition, but satisfies
the bound $\rho_t^{(n)} (x) \leq C_n t^{-n/2}$ uniformly over all possible finite initial conditions $(x_i)$. This follows
by duality for $n=1$ and by anti-correlation for $n>1$ (see \cite{mrtz}). It follows that
$E^{C}_{(x_i)} [N_t^p(a,b)]$ is bounded, for each $t,p>0,\, a,b \in \R$, uniformly over finite initial conditions $(x_i)$.

Fix $\mu \in \mathcal{M}_0$ and take finite measures $\mu_n$ so that $\mu_n \to \mu$ (recall we are using vague convergence).
The moment bounds above imply that the laws of $N_t$ on $\mathcal{M}_{LFP}$ under $P^C_{\mu_n}$ are tight. Take a subsequence $n'$ along which
they converge to a limit, which we denote $Q$. The functions
\[
\nu \to F_{(a_i)} (\nu) :=  \chi(\nu(I_1) = \nu(I_3) = \ldots = \nu(I_{2n-1})=0)
\]
are discontinuous on $\mathcal{M}_{LFP}$. However the moment bound $E^{C}_{(x_i)} [N_t(a,b)] \leq C(t) (b-a)$
holds also for the limit law $Q$ and implies that
$\nu(\{a_i\})=0$, $Q(d\nu)$ almost surely. This shows that $Q$ does not charge the discontinuity set of $F_{(a_i)}$.
Then we may pass to the limit in  (\ref{cbmduality}) to deduce that
\begin{equation} \label{info}
\int F_{(a_i)}(\nu) Q(d \nu) = P^{A}_{(a_i)} \left[ S_t \cap \mbox{supp}(\mu) = \emptyset \right].
\end{equation}
These functionals do not characterize a law on $\mathcal{M}_{LFP}$, but they do characterize a law that is
supported on $\mathcal{M}_0$. To see this note that for $\nu \in \mathcal{M}_0$
\[
\nu([x,y]) = \lim_{N \to \infty} \sum_k \chi \left( \nu ([x,y] \cap (\small{\frac{k}{N}, \frac{k+1}{N}}]) >0 \right).
\]
From this one may use (\ref{info}) to
find $\int \nu([x_1,y_1]) \ldots \nu([x_n,y_n]) Q(d \nu)$ which, by the moment bounds, determine $Q$.
To see that $Q$ is supported on $\mathcal{M}_0$ note that
\[
P^C_{\mu_n} \left[ N_t(a,b) \geq 2 \right]  \leq \int_a^b \int_a^b \rho^{(2)}_t(x_1,x_2) \, dx_1 dx_2
\leq C(t) (b-a)^2.
\]
This bound holds uniformly over $n$ and hence also for the limit law $Q$. Then the conclusion follows
from the usual covering argument, for instance
\begin{eqnarray*}
&& \hspace{-.3in} Q \left[ \mbox{ $\mu(\{x\} >1$ for some $x \in [-L,L]$)} \right] \\
& \leq &  \sum_{k=-LN}^{LN}
Q \left[ \mu( [k/N,(k+1)/N]) \geq 2 \right] \leq C(L,t) N^{-1}.
\end{eqnarray*}
Thus the law $Q$ is determined and we may define $p_t(\mu,d\nu)$ to equal $Q(d\nu)$.

The remainder of the results in section \ref{s2.3} follow using similar tools. For example, for the
continuity of $\mu \to p_t(\mu,d\nu)$, that is the Feller property, suppose that $\mu_n \to \mu$
in $\mathcal{M}_0$. The moment bounds, which still hold for infinite initial conditions,
imply the tightness of $p_t(\mu_n,d\nu)$. Passing to the limit in
\[
\int F_{(a_i)}(\nu) p_t(\mu_n, d\nu) = P^{A}_{(a_i)} \left[ S_t \cap \mbox{supp}(\mu_n) = \emptyset \right].
\]
shows that any limit point of $p_t(\mu_n,d\nu)$ must be $p_t(\mu,d\nu)$. The semigroup property, for bounded
 continuous $F: \mathcal{M}_{LFP} \to \R$,
\begin{equation} \label{semigroup}
\int F(\nu) p_{t+s}(\mu,d\nu) = \int \int F(\nu') p_s(\nu,d\nu') p_t(\mu,d\nu),
\end{equation}
which is valid for finite measures $\mu$ extends  to hold for $\mu \in \mathcal{M}_0$ by approximation, using the Feller
property.
The same tightness and characterization methods establish the existence of a law
characterized by (\ref{cbmentrance}), and justify the arguments in (\ref{attraction})
that many initial laws are attracted to it. That (\ref{cbmentrance}) determines an entrance
law can be established by passing to the
limit in (\ref{semigroup}) along $\mu = \sum_k \delta_{\lambda^{-1}k}$ as $\lambda \to \infty$.

The annihilating case follows the same lines, with moments controlled since the $n$-point density and moments
for annihilating systems are bounded by the corresponding coalescing system. The coalescing duality formula (\ref{cbmduality})
is replaced by the annihilating duality formula (\ref{abmduality}), and to see that this will characterize the law note that
for $\nu \in \mathcal{M}_0$
\[
\nu([x,y]) = \lim_{N \to \infty} \sum_k \left( 1 - (-1)^{\nu \left([x,y] \cap (\frac{k}{N},\frac{k+1}{N}] \right)} \right).
\]
\subsection{Details for Remark 3 in section \ref{s3.3}} \label{s4.2}
We give here the algebraic manipulations to derive the Pfaffians (\ref{cbmpf2},\ref{cbmpf2.5},\ref{cbmpf3},\ref{cbmpf4}).

$\mbox{Pf}\,(I)=1$ for $I$ the $2n \times 2n$ anti-symmetric matrix with entries $1$ above the diagonal,
and the formula (\ref{pfaffsum}) specializes to
\[
\mbox{Pf}\,(I-A) = \sum_{J} (-1)^{s(J)} \mbox{Pf}\,(A|_{J}),
\]
(using for a $2n \times 2n$ anti-symmetric matrix $A$, that $\mbox{Pf}\,(-A) = (-1)^n \mbox{Pf}\,(A)$).
We combine this with a simple combinatorial identity (which can be checked
 by induction on $n$): suppose that $(m_{j,k}:1 \leq j < k \leq n)$ satisfy the collapsing product
$m_{j,k}m_{k,l}=m_{j,l}$ for all $j,k,l$; then
\begin{eqnarray*}
\prod_{k=1}^{n-1} (1+m_{k,k+1}) & = & 1 + \sum_{1 \leq k_1 < k_2 \leq n} m_{k_1,k_2}  \\
&& \hspace{.2in} + \sum_{1 \leq k_1 < k_2 <k_3 <k_4 \leq n} m_{k_1,k_2} m_{k_3,k_4}  \\
&& \hspace{.2in} \hspace{.2in} + \sum_{1 \leq k_1 < k_2 <k_3 <k_4 <k_5<k_6 \leq n} m_{k_1,k_2} m_{k_3,k_4} m_{k_5,k_6} + \ldots \\
& = & \sum_{J} m_J
\end{eqnarray*}
where the final sum is over all subsets $J$ of $\{1,2,\ldots,n\}$ of even size,
and if $J=\{k_1,\ldots,k_{2m}\}$ where $k_1< \ldots <k_{2m}$ then
$m_J = m_{k_1,k_2} m_{k_3,k_4}  \ldots m_{k_{2m-1},k_{2m}}$ (and with $m_{\emptyset}=1$).
If $n$ is even then the last term of this series is $m_{1,2}m_{3,4} \ldots m_{n-1,n}$.
 Note that $\bar{m}_{j,k} = \alpha^{k-j}m_{j,k}$ also satisfy
$\bar{m}_{j,k} \bar{m}_{k,l} = \bar{m}_{j,l}$ and applying the above for $\bar{m}$ one obtains a
decomposition for $\prod_{k=1}^{n-1} (1+ \alpha m_{k,k+1})$. In particular for $\alpha=-1$ we get
\[
\prod_{k=1}^{N-1} (1-m_{k,k+1}) = \sum_{J} (-1)^{s(J)} m_J.
\]
Now apply this with $m_{j,k} = \chi(N_t \left((a_j,a_k)\right)=0)$. These satisfy the collapsing products almost surely under the
probability $P^{C}_{\infty}$.  The Pfaffian (\ref{cbmpf}) shows that
$E^{C}_{\infty} [m_J] =  \mbox{Pf}\,(F|_J)$ and so
\begin{eqnarray}
P^{C}_{\infty} \left [\mbox{$N_t(I_k) > 0 $ for $k=1,\ldots,2n-1$} \right]
& = & E^{C}_{\infty} \left [\prod_{k=1}^{2n-1} (1-m_{k,k+1}) \right] \nonumber \\
& = &  \sum_{J} (-1)^{s(J)} E^{C}_{\infty} \left [m_J \right] \nonumber \\
& = &  \sum_{J} (-1)^{s(J)} \mbox{Pf}\,(F|_J) \nonumber \\
& = & \mbox{Pf}\,(I-F). \nonumber
\end{eqnarray}
We may apply the same argument for the annihilating case taking $m_{j,k} = (-1)^{N_t \left((a_j,a_k)\right)}$,
where $1-m_{j,k} = 2 \chi(\mbox{$N_t \left((a_j,a_k)\right)$ is odd})$, to find (\ref{cbmpf2.5}).

For (\ref{cbmpf3}) we have $\mbox{Pf}\,(O_{2n}) = 1$ and $\mbox{Pf}\,(O_{2n}|_{J}) = 0$ unless $O_{2n}|_{J}$ is a copy of $O_{2m}$ for
some $m \in \{0,1,\ldots,n\}$. This occurs either if $J$ is empty or if $J$ is of the form
\begin{eqnarray}
J_1 & = & \{2k_1-1,2k_1, 2k_2-1,2k_2,\ldots,2k_m-1, 2k_m\} \nonumber \\
& & \quad \mbox{for some $1 \leq k_1 < \ldots <k_m \leq n$}.
\label{specialform}
\end{eqnarray}
Then formula (\ref{pfaffsum}) specializes to
\[
\mbox{Pf}\,(O-A) = \sum_{J_1} (-1)^{|J_1|/2} \mbox{Pf}\,(A|_{J_1})
\]
where the sum is over all $J_1$ of the form in (\ref{specialform}) (including the empty set). We use another
combinatorial identity, also straightforward by induction on $n$:
\[
\prod_{k=1}^{n} \left(1-m_{2k-1,2k} \right) = \sum_{J_1} (-1)^{|J_1|/2} m_{J_1}
\]
where the sum is over all $J_1$ of the form in (\ref{specialform}) (including the empty set). Arguing as in the
previous example leads to  (\ref{cbmpf3}).

For (\ref{cbmpf4}) one has $\mbox{Pf}\,(\hat{O}_{2n}) = 0$ and $\mbox{Pf}\,(\hat{O}_{2n}|_{J}) = 0$ unless $\hat{O}_{2n}|_{J}$ is a copy of
$O_{2m}$ for some $m \in \{0,1,\ldots,n-1\}$. This occurs either if $J$ is empty or if $J$ is of the form
\begin{eqnarray}
J_2 & = & \{2k_1,2k_1+1, 2k_2,2k_2+1,\ldots,2k_m, 2k_m+1\} \nonumber \\
& & \quad \mbox{for some $1 \leq k_1 < \ldots <k_m \leq n-1$}.  \label{specialform2}
\end{eqnarray}
Then formula (\ref{pfaffsum}) specializes to
\[
\mbox{Pf}\,(A-\hat{O}) = \sum_{J_2} (-1)^{|J_2|/2} \mbox{Pf}\,(A|_{J_2}^c)
\]
where the sum is over all $J_2$ of the form in (\ref{specialform2})  (including the empty set).
The required combinatorial identity is
\[
\prod_{k=1}^{n-1} \left(1-m_{2k,2k+1} \right) \, \prod_{k=1}^n m_{2k-1,2k} = \sum_{J_2} (-1)^{|J_2|/2} m_{J_2}
\]
where the sum is over all $J_2$ of the form in (\ref{specialform2})  (including the empty set).
Arguing as in the previous examples leads to  (\ref{cbmpf4}).
\subsection{Details on distributional derivatives} \label{s5.3} \label{s4.3}
The derivation of the kernel $K$ in section \ref{s3.3}, and also the Pfaffian (\ref{tildeclaim}), use formal differentiation
that can be made precise by using distributional derivatives. Consider first  (\ref{tildeclaim}).
For  $\mu = \sum_i \delta_{z_i}$ a locally finite point measure with disjoint atoms, one has
 the distributional derivative on $V_{2n}$
\begin{eqnarray}
&& \hspace{-.3in} \partial_{x_1} \ldots \partial_{x_{2n}}
\chi(\mbox{$\mu(x_k,x_{k+1})=0$ for $k=1,2,\ldots,2n-1$}) \nonumber \\
&=& (-1)^n \chi(\mbox{$\mu(x_k,x_{k+1})=0$ for $k=1,3,5\ldots,2n-1$}) \mu(dx_1) \ldots \mu(dx_{2n}). \label{1001}
\end{eqnarray}
We illustrate how to check this by showing that, in the distributional sense on $V_2$,
\[
\partial_x \chi(\mu(x,y)=0) = \chi(\mu(x,y)=0) \mu(dx)  \, dy.
\]
Indeed, if $f$ is smooth and compactly supported in $V_2$ then
\begin{eqnarray*}
&& \hspace{-.3in}
 \int_{\R^2} f(x,y) \chi(\mu(x,y)=0) \mu(dx)  \, dy \\
& = & \sum_i \int_{\R} f(x_i,y) \chi(\mu(z_i,y)=0, \, z_i< y ) \, dy \\
& = & \int_{\R^2} \partial_x f(x,y) \left( \sum_i  \chi(\mu(z_i,y)=0,\, x < z_i < y) \right) \,dx \,  dy \\
& = & \int_{\R^2} \partial_x f(x,y) \chi(\mu(x,y)>0) \, dx  \, dy \\
& = & - \int_{\R^2} \partial_x f(x,y) \chi(\mu(x,y)=0) \, dx  \, dy
\end{eqnarray*}
since at most one term in the sum over $i$ is non-zero. Iterating such calculations leads to  (\ref{1001}).
Then for smooth $f$ compactly supported in $V_{2n}$,
\begin{eqnarray*}
&&  \hspace{-.4in} \int_{V_{2n}} f(x_1,\ldots,x_{2n}) \tilde{\rho}^{(2n)}_t(x_1,\ldots,x_{2n}) dx_1 \ldots dx_{2n} \\
& = &  E^{C}_{\infty} \left[ \int_{V_{2n}}  f(x_1,\ldots,x_{2n})
 \chi(\mbox{$N_t(I_k)=0$ for $k=1,3,\ldots,2n-1$}) N_t(dx_1) \ldots N_t(dx_{2n})  \right] \\
& = & (-1)^n E^{C}_{\infty} \left[ \int_{V_{2n}}  \partial_{x_1} \ldots \partial_{x_{2n}} f(x_1,\ldots,x_{2n})
 \chi(\mbox{$N_t(I_k)=0$ for $k=1,2,\ldots,2n-1$}) dx_1 \ldots dx_{2n}  \right] \\
& = &  (-1)^n \int_{V_{2n}}  \partial_{x_1} \ldots \partial_{x_{2n}} f(x_1,\ldots,x_{2n})
\mbox{Pf}\, \left(F(t^{-1/2}(x_j-x_i) )\right) \, dx_1 \ldots dx_{2n}  \\
& = &  (4 \pi t^2)^{-n/2} \int_{V_{2n}}  f(x_1,\ldots,x_{2n})
\mbox{Pf}\,\left(\phi(t^{-1/2}(x_j-x_i))\right) \, dx_1 \ldots dx_{2n}.
\end{eqnarray*}
In the last step we have passed the derivatives onto the Pfaffian, which is smooth since $F$ is smooth, and used
$F''(x) = (4 \pi)^{-1/2} \phi(x)$.

The argument for the kernel $K$ is similar. Fix $x_2<x_4< \ldots <x_{2n}$
and consider the open set $V = \{ (x_1,x_3,\ldots,x_{2n-1}): x_1 < x_3 < \ldots < x_{2n-1}\}$.
Then, as above, in the distributional sense on $V$
\begin{eqnarray}
&& \hspace{-.3in} \partial_{x_1} \partial_{x_3} \ldots \partial_{x_{2n-1}}
\chi(\mbox{$\mu(x_k,x_{k+1})=0$ for $k=1,3,\ldots,2n-1$}) \nonumber \\
&=& \chi(\mbox{$\mu(x_k,x_{k+1})=0$ for $k=1,3\ldots,2n-1$}) \mu(dx_1) \, \mu(dx_3) \ldots \mu(dx_{2n-1}). \nonumber
\end{eqnarray}
Then for smooth $f$ compactly supported in $V$, with
$\Omega = \{\mbox{$N_t(I_k)=0$ for $k=1,3,\ldots,2n-1$}\}$,
\begin{eqnarray*}
&&  \hspace{-.4in}  E^{C}_{\infty} \left[ \int_{V} f(x_1,x_3,\ldots,x_{2n-1})
 \chi(\Omega) N_t(dx_1) \, N_t(dx_3)  \ldots N_t(dx_{2n-1})  \right] \\
& = & (-1)^n E^{C}_{\infty} \left[ \int_{V}  \partial_{x_1} \partial_{x_3} \ldots \partial_{x_{2n-1}}
f(x_1,x_3, \ldots,x_{2n-1})
 \chi(\Omega) dx_1 \, dx_3 \ldots dx_{2n-1}  \right] \\
& = &  (-1)^n \int_{V}   \partial_{x_1} \partial_{x_3} \ldots \partial_{x_{2n-1}}
f(x_1,x_3, \ldots,x_{2n-1})     \\
&& \hspace{.5in}
\mbox{Pf}\,\left(F(t^{-1/2}(x_j-x_i))\,:\, 1 \leq i < j \leq 2n \right) \, dx_1 \, dx_3 \ldots dx_{2n-1}  \\
\end{eqnarray*}
 Now one can pass the derivatives onto the Pfaffian and then let $x_{2} \downarrow x_{1}, \, x_4 \downarrow x_3, \ldots $
as described in section \ref{s3.3}.
\subsection{Details for section \ref{s3.4}} \label{s4.4}
Here we give the error estimates for the Pfaffian expansion Lemma \ref{olegslemma}.

The product (\ref{symmetric}) that defines the symmetric polynomials $\sigma^{2n}_k$
yields a total of $2^{2n}$ monomials so we have the simple bound $|\sigma_k^{2n}(y)| \leq 2^{2n} |y|^k$.
The expansion (\ref{pretau}) must be iterated at most $p$ times to derive (\ref{tau}) and this leads to
to the bound
\begin{equation} \label{taubound}
|\tau^{2n,p+2n}_k(y)| \leq (2n 2^{2n})^p |y|^{p+2n-k+1}.
\end{equation}
Using this we may bound the size of the remainder terms given in (\ref{Reqn}). For example
\begin{eqnarray*}
&& \hspace{-.3in}  \sum_{k,l=2n+1}^{\infty} | \tau^{2n,l-1}_p(y)| \, | J_{kl} | \, | \tau^{2n,k-1}_q(y) | \\
& \leq &  \sum_{k,l=2n+1}^{\infty} \frac{1}{(k-1)! (l-1)!} \, |\phi^{k+l-2}(0)| \, (2n 2^{2n})^{l+k-4n-2} |y|^{l+k-p-q} \\
& \leq & |y|^2 \sum_{k,l=2n+1}^{\infty} \frac{1}{(k-1)! (l-1)!} \, |\phi^{k+l-2}(0)| \, (2n 2^{2n} \epsilon)^{l+k-4n-2} \quad
\mbox{when $|y| \leq \epsilon$} \\
& \leq & |y|^2 \sum_{r=4n}^{\infty} \sum_{|s| \leq r-4n}
\frac{2^r}{r!} \, |\phi^{r}(0)| \, (2n 2^{2n} \epsilon)^{r-4n} \\
&& \qquad \mbox{using $r=k+l-2$, $s=k-l$ and $\frac{k! \, l!}{(k+l)!} \geq 2^{-k-l}$} \\
& \leq & 2^{4n} |y|^2  \sum_{r=4n}^{\infty}
\frac{1}{r!} \, |\phi^{r}(0)| \, 2r (4n 2^{2n} \epsilon)^{r-4n}.
\end{eqnarray*}
Choosing $\epsilon = \epsilon(n,\phi)$ so that $4n 2^{2n} \epsilon$ lies in the radius of convergence of $\phi$
we obtain a convergent series. Similarly
\begin{eqnarray*}
 \sum_{l=2n+1}^{\infty} | \tau^{2n,l-1}_p(y)| \, | J_{ql} |
& \leq &   \sum_{l=2n+1}^{\infty} \frac{1}{(q-1)! \,(l-1)!} \, |\phi^{q+l-2}(0) | \, (2n2^{2n})^{l-2n-1} |y|^{l-p} \\
& \leq &  |y|  \sum_{l=2n+1}^{\infty} \frac{1}{(l-1)!} \, |\phi^{q+l-2}(0) | \, (2n2^{2n} \epsilon)^{l-2n-1} \\
& \leq &  |y|  \sum_{l=2n+1}^{\infty} \frac{1}{(q+l-2)!} \, |\phi^{q+l-2}(0) | \, (2n2^{2n} \epsilon)^{l-2n-1} (l+2n)^{2n} \\
& \leq & C(n,\phi) |y|.
\end{eqnarray*}
A similar bound holds for the final term in (\ref{Reqn}). Combining the estimates yields the desired error bound on $R^{(2n)}_{pq}$.
Moreover these bounds show the absolute convergence that justifies the rearrangement of the series (\ref{newton})
 used in Lemma \ref{olegslemma} provided that
$|y| \leq \epsilon(n,\phi)$.

\end{document}